# Optimization of Survivability Analysis for Large-Scale Engineering Networks

S.V. Poroseva

Department of Mechanical Engineering, University of New Mexico,
Albuquerque, NM 87131-0001, U.S.A.

P.A. Rikvold

Department of Physics, Florida State University,
Tallahassee, FL 32306-4350, U.S.A.

**Abstract**

Engineering networks fall into the category of large-scale networks with heterogeneous nodes such as sources and sinks. The survivability analysis of such networks requires the analysis of the connectivity of the network components for every possible combination of faults to determine a network response to each combination of faults. From the computational complexity point of view, the problem belongs to the class of exponential time problems at least. Partially, the problem complexity can be reduced by mapping the initial topology of a complex large-scale network with multiple sources and multiple sinks onto a set of smaller sub-topologies with multiple sources and a single sink connected to the network of sources by a single link. In this paper, the mapping procedure is applied to the Florida power grid.

**Keywords:** networks, survivability, resilience, computational complexity

## 1 Introduction

Traditionally, engineering networks such as, for example, electric power, gas, water, transportation systems etc. are expected to operate under normal conditions most of the time. Their behavior in the presence of predictable faults with the predictable time of fault repair can be described using reliability/availability analysis.

In contrast, the ability of a network to withstand massive sudden damage caused by adverse events (or *survivability*) has not been among design goals in the past. In the modern society, this practice needs an adjustment. Indeed, modern networks develop toward increasing their size, complexity, and integration. The likelihood of adverse events increases as well due to technological development, climate change, and activities in the political arena among other factors. Under such circumstances, a network failure has an unprecedented effect on lives and economy. Some networks, for example, telecommunication and financial systems, have already reached the global level. World economical crises are an unfortunate example of how failures in hyper-large networks can affect societies worldwide.

To mitigate the consequences of adverse events, mathematical and computational tools are required to evaluate the survivability of existing networks, analyze the efficiency of different design strategies for enhancing the network survivability, and compare the survivability of alternative designs for future networks. The development of such tools is a goal of our research.

Many factors contribute to the network's survivability. Our study analyzes the impact of the network's topology on its survivability. The network's topology is defined as the number of network components, their type, and how the network components are connected to one another.

Engineering networks contain many different components, but for the survivability analysis, only four types of components are of importance: three types of nodes and links that connect the nodes with one another. The three types of nodes are: sources that generate a quality or service of interest, sinks that consume this quality or service, and interconnections through which the quality/service is transported without a change (ideally) in its amount in the direction of a flow of the quality/service. The flow direction through interconnections may vary depending on the network state.

In such networks, survivability is perceived as the continuation of a supply of a quality of interest (service) from sources to sinks in an amount sufficient to satisfy their demand in the presence of multiple faults in network elements. Depending on the origin and evolution of faults, the amount of a quality of interest available to each sink may vary with time. Intuitively, however, the word "survivability" is associated not with a process, but with a final state of a system after all possible damage has occurred and before any recovery action can take place. This is the approach adopted in our study, that is, only a final steady state of a network is considered after all faults (including cascaded and secondary faults) have occurred and before any repair has been accomplished. As the origin and the development of faults are not of concern, the research is applicable to adverse events of any type and intensity in which multiple faults have resulted either due to multiple events (independent or correlated), or due to a single event causing cascaded/secondary faults in multiple network elements.

In the final steady state, one determines whether a network survived an adverse event by comparing the amount of the quality of interest (service) available in the transformed network with the demand for this quality/service. A mathematical problem to be solved is to identify sources that survived faults, determine their capacities, and verify their connections to sinks. The problem has to be solved for all possible combinations of faults (or *fault scenarios*).

Previously, we suggested a probabilistic framework for quantifying the network survivability [1] and developed an efficient numerical algorithm for conducting the survivability analysis of small- to medium-sized networks with multiple sources and a single sink [2]. This approach is applicable, for example, to power grids, when a single sink represents either an isolated industrial load, or multiple commercial and residential loads interconnected into a single distribution system, or a lower-voltage level network.

In larger and more complex networks, the algorithm becomes computationally unfeasible. The main factor is the size of the problem. Indeed, in a network with $M$ elements, the total number of possible final steady states is $2^M$. (In a final steady state, each element can only be in one of the two states: available or faulty.) From the computational complexity point of view, the problem belongs to the class of exponential time problems at least [3,4]. That is, the problem is more complex than the more familiar NP-complete problems.

Some ways to reduce the problem complexity were suggested in [2,5]. In [6], a "selfish" algorithm was described, in which the problem complexity is reduced by mapping the initial topology of a complex large-scale network with multiple sources and sinks onto a set of simpler smaller topologies (or *sub-topologies*) with multiple sources and a single sink connected to the grid by a single link.

In [6], the algorithm was used to disintegrate the topology of the notional Medium Voltage DC shipboard power system [7] into a set of sub-topologies. It was shown that after all simplifications, the initial topology to analyze contains 32 network components that represent four generators, seven loads, and links. That is, the initial topology is relatively small. Yet, the total number of possible faults scenarios to generate and search for the components connectivity is $2^M$ or about 4 billion. By applying the "selfish" algorithm, however, the initial topology can be mapped to a single sub-topology that contains only 10 components. In this sub-topology, only 1024 fault scenarios have to be generated and searched.

Whereas the "selfish" algorithm does not eliminate the need to generate all fault scenarios in the initial topology, it substitutes a search procedure applied to every fault scenario in this network by a decision problem. In the decision problem, the network response to a given combination of faults is determined by requesting information from a database previously created for a set of simpler smaller sub-topologies (see [4] for more discussion on search and decision problems).

The current paper describes our initial results in mapping the Florida power grid onto a set of sub-topologies.

## 2  Network representation

In our previous works [5,6], it was found that for the survivability analysis of networks with sources and sinks, the most efficient way of representing a network is by links only. Indeed, if a source or a sink is connected to the network by a single link, the node and the link are connected in series. Then, a fault in the link makes the node unavailable to the network. A fault in the node makes the link useless. Therefore, the node and the link can be represented by a single link.

If there are several links connecting a source or a sink to the network, faults in all of them isolate the node from the network, which is equivalent to a fault in the node. Thus, there is no need to consider faults in the nodes separately, and they can be removed from the network. Similarly, faults in interconnections can also be removed from consideration.

In the network representation by links, there are three different types of links: links adjacent to sources (vertical "VT" links), links adjacent to sinks (vertical "VB" links), and links between interconnections (horizontal "H" links). To recognize the three types of links mathematically, the VB- and VT-links are assigned weights. The sign of the weight shows the flow direction through a link: the VT-links have positive weight and the VB-links have negative weight. No sign is assigned to the H-links, because the flow direction through such links may vary depending on the network topology after faults. A value of the weight assigned to a vertical link can be based, for example, on its capacity or a maximum amount of the quality of interest (service) it can transport in relation to the total demand.

To demonstrate the network representation by links only, let us consider as an example a small part of the Florida power grid. In Figure 1a, the circles with "-" correspond to sinks (loads), the circles with "+" show sources (power plants, generators), and the dark circles are distribution substations. The generators and distribution substations are numbered in accordance with the Florida power grid map (Fig. 4).

Usually, power system diagrams do not explicitly show loads, but only links to them from distribution substations. The map in Fig. 4 also includes only generators and distribution substations. In Figure 1a, however, loads are shown (but not numbered) to bring more clarity in the explanation of the network representation procedure.

In the network in Fig. 1a, elements connected in series can be represented as a single element. For example, a single VB-link (VB30) can be used to represent node 30, links that are adjacent to the node, and the load next to this node (Fig. 1b). The other loads with the links adjacent to them can also be substituted by VB-links: VB20, VB27, and VB28. The number of a VB-link corresponds to that of the distribution substation to which the load is adjacent.

Nodes 64, 76, and 81 that represent generators in Fig. 1a and the links adjacent to them can be substituted with VT-links. For nodes 64 and 81, these are links VT64 and VB81 in Fig. 1b. The two links adjacent to node 76 are shown as VT761 and VT762. Nodes 20, 27, 28, and 30 (distribution substations) become interconnections (not numbered) in Fig. 1b. Links between the interconnections are H-links: H1 and H2. In the figure, the arrow heads of links show the flow direction.

The network representation shown in Fig. 1b is the one used in the survivability analysis.

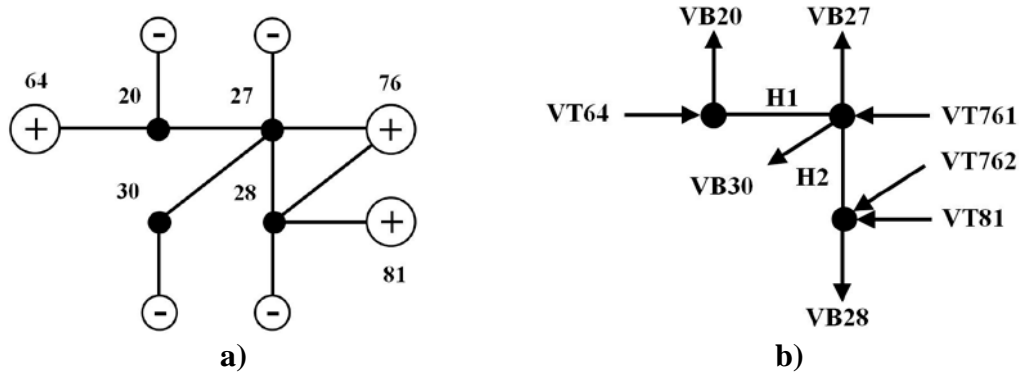

Figure 1: a) initial network representation; b) network representation for the survivability analysis.

## 3 Mapping procedure

A general procedure for disintegrating a network with multiple sources and sinks into a set of sub-topologies is shown in Fig. 2. The network disintegration is achieved by first removing all VB-links belonging to other sinks from the initial network. Indeed, faults in vertical links connecting other sinks to the network cannot interrupt the flow of the quality of interest (service) to the sink under consideration. Then, in each sub-topology, links connected in series are combined into a single link. Since faults in links connected in series are equivalent to a fault in a single link, this procedure is justified for the survivability analysis. Moreover, it reduces the number of links in a sub-topology.

In a general case, the number of sub-topologies is equal to the number of sinks in the initial network. However, one can expect that many sinks "see" the network alike, and therefore the number of sub-topologies to consider will be much less.

The last step in this procedure is the disintegration of sub-topologies with a sink connected to the network by multiple VB-links into sub-topologies with a sink connected to the network by a single VB-link. Not only does this step reduce the scale of the sub-topologies further, it will also, most likely, reduce their number due to similarity between many of them.

In the worst-case scenario, when no similarity between sub-topologies is found, the final number of sub-topologies to consider is equal to the number of VB-links, that is, *VB*. As only one VB-link is included in a sub-topology, each sub-topology can be uniquely identified by the VB-link.

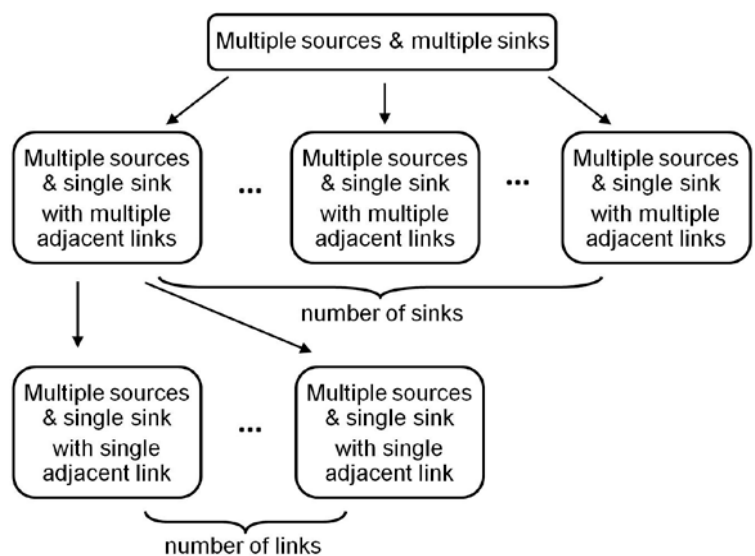

Figure 2: Mapping procedure.

Let us consider as an example the disintegration of the network shown in Fig. 1b. The initial topology with multiple sources and multiple sinks contains 10 elements: four VB-links, four VT-links, and two H-links.

A set of sub-topologies of this network should only include networks with multiple VT-links and a single VB-link. Figure 3a shows a sub-topology for link VB20 obtained by removing other VB-links from the initial network. This sub-topology contains seven elements.

The sub-topology for link VB27 is shown in Fig. 3b. This sub-topology can further be simplified by representing two links in series – VT64 and H1 – as a single link VT64 (Fig. 3c). The sub-topology in Fig. 3c is the final one for link VB27 to analyze. This sub-topology contains six elements.

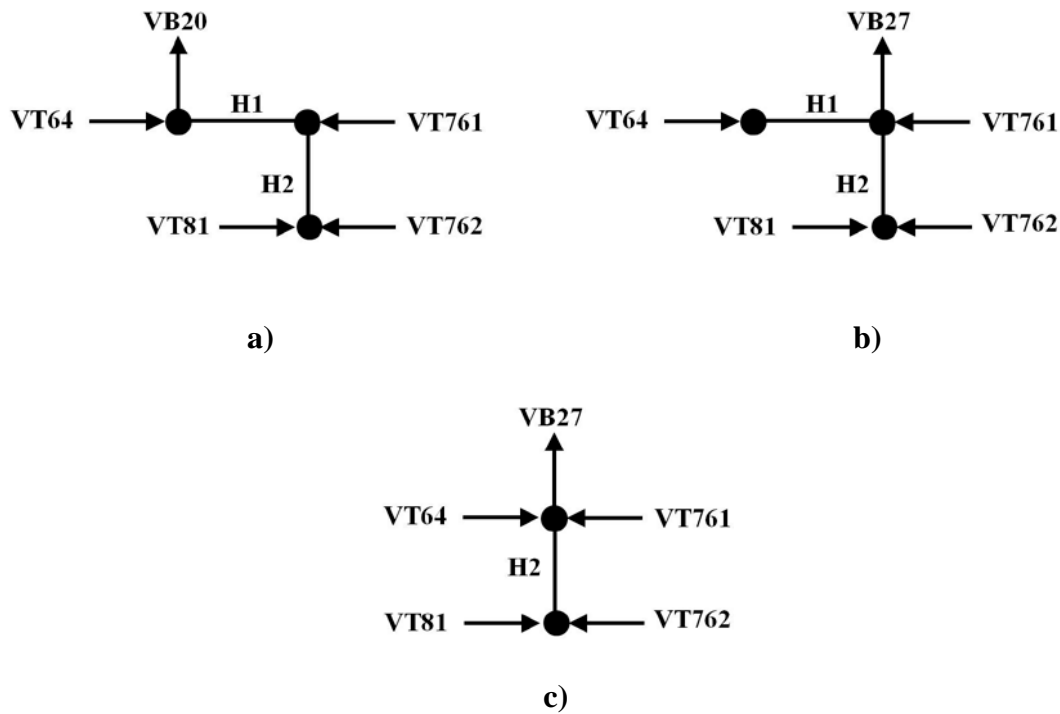

Figure 3: Sub-topologies for links VB20 and VB27 in the network in Fig. 1b.

Similar steps conducted for links VB28 and VB30 result in the same sub-topology as the one in Fig. 3c, where link VB27 is substituted with link VB28 or link VB30.

To conclude, the initial topology in Fig. 1b with $2^{10} = 1024$ fault scenarios can be mapped onto a set of two sub-topologies (Figs. 3a and 3c). The combined number of fault scenarios in these two sub-topologies is $2^7 + 2^6 = 192$, that is, the number of fault scenarios to analyze is approximately five times less than in the initial network.

## 4  Application to the Florida power grid

The Florida high-voltage power grid [8] contains $M = 84$ nodes, 31 of which are generators and the rest are distribution substations that act as loads. The nodes are connected by 200 links, some of which are parallel power lines connecting the same two nodes.

Figure 4 shows a map of the Florida grid [9,10]. In the figure, squares correspond to generators and circles to loads. For the purposes of the survivability analysis, if two nodes are connected by parallel links that are geographically co-located, such links are represented as a single link. Indeed, external damage due to adverse events will most likely occur simultaneously in all co-located links. Therefore, the map does not contain parallel links.

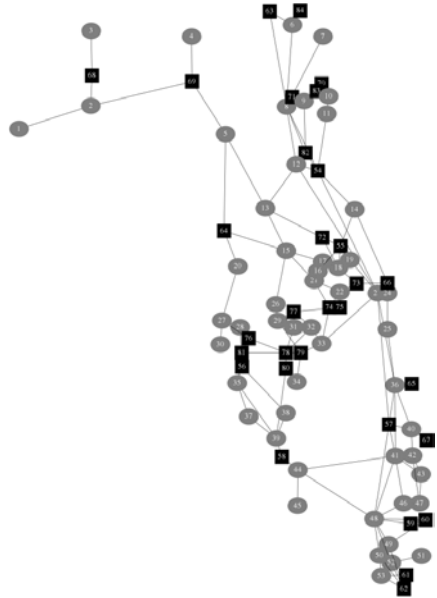

Figure 4: Florida power grid

As one can see, the grid's topology is quite complex, even without parallel links. A reason is that in contrast to the shipboard power system [7] and other microgrids, utility systems such as the Florida power grid are not designed as a whole system, but rather evolve following economic development in a region. The analysis of such a system is also complicated due to the absence of obvious symmetries in the network topology, insufficient and incomplete data about the system, the large number of elements and their complicated connections with one another. In this regard, the current paper is the first step in the analysis of such a grid.

The preliminary analysis revealed that the system can be disintegrated in a few groups of interconnected generators and distribution substations.

| **Group 1** | **Group 2** | **Group 3** | **Group 4** | **Group 5** |
|---|---|---|---|---|
| 3:{68} | 4:{69} | 34:{79,80} | 14:{54,55,66} | 9:{70,71,82,83} |

| **Group 6** | **Group 7** | **Group 8** | **Group 9** |
|---|---|---|---|
| 2:{68,69} | 10:{54,70,83} | 27: {64,76,81} | 35:{56,58,80,81} |
| 1:{68,69} | 11:{54,70,83} | 30:{64,76,81} | 37: {56,58,80,81} |
| | | 20: {64,76,81} | 38: {56,58,80,81} |
| | | 28: {64,76,81} | 39:{56,58,80,81} |

**Group 10**
15: {55,64,69,72,73,74,75,77,79,82]
21: {55,64,69,72,73,74,75,77,79,82}
17:{ 55,64,69,72,73,74,75,77,79,82}
16: {55,64,69,72,73,74,75,77,79,82}
26: {55,64,69,72,73,74,75,77,79,82}
13: {55,64,69,72,73,74,75,77,79,82)
5:{ 55,64,69,72,73,74,75,77,79,82}
29: {55,64,69,72,73,74,75,77,79,82}
18:{ 55,64,69,72,73,74,75,77,79,82}
22: {55,64,69,72,73,74,75,77,79,82}

**Group 11**
40:{ 54,55,57,58,59,60,61,62,63,65,66,67,71,74,77,78,79,82,84}
41: {54,55,57,58,59,60,61,62,63,65,66,67,71,74,77,78,79,82,84}
48:{ 54,55,57,58,59,60,61,62,63,65,66,67,71,74,77,78,79,82,84}
49:{ 54,55,57,58,59,60,61,62,63,65,66,67,71,74,77,78,79,82,84}
50:{ 54,55,57,58,59,60,61,62,63,65,66,67,71,74,77,78,79,82,84}
36:{ 54,55,57,58,59,60,61,62,63,65,66,67,71,74,77,78,79,82,84}
23:{ 54,55,57,58,59,60,61,62,63,65,66,67,71,74,77,78,79,82,84}
8:{ 54,55,57,58,59,60,61,62,63,65,66,67,71,74,77,78,79,82,84}
12:{ 54,55,57,58,59,60,61,62,63,65,66,67,71,74,77,78,79,82,84}
46:{ 54,55,57,58,59,60,61,62,63,65,66,67,71,74,77,78,79,82,84}
42:{ 54,55,57,58,59,60,61,62,63,65,66,67,71,74,77,78,79,82,84}
47:{ 54,55,57,58,59,60,61,62,63,65,66,67,71,74,77,78,79,82,84}
43:{ 54,55,57,58,59,60,61,62,63,65,66,67,71,74,77,78,79,82,84}
24:{ 54,55,57,58,59,60,61,62,63,65,66,67,71,74,77,78,79,82,84}
44:{ 54,55,57,58,59,60,61,62,63,65,66,67,71,74,77,78,79,82,84}
52:{ 54,55,57,58,59,60,61,62,63,65,66,67,71,74,77,78,79,82,84}
25:{ 54,55,57,58,59,60,61,62,63,65,66,67,71,74,77,78,79,82,84}
33:{ 54,55,57,58,59,60,61,62,63,65,66,67,71,74,77,78,79,82,84}
32:{ 54,55,57,58,59,60,61,62,63,65,66,67,71,74,77,78,79,82,84}
6:{ 54,55,57,58,59,60,61,62,63,65,66,67,71,74,77,78,79,82,84}
31:{ 54,55,57,58,59,60,61,62,63,65,66,67,71,74,77,78,79,82,84}
19:{ 54,55,57,58,59,60,61,62,63,65,66,67,71,74,77,78,79,82,84}
7:{ 54,55,57,58,59,60,61,62,63,65,66,67,71,74,77,78,79,82,84}
45:{ 54,55,57,58,59,60,61,62,63,65,66,67,71,74,77,78,79,82,84}

Here, nodes from 1 to 53 are distribution substations and nodes from 54 to 84 are generators. The connectivity of the network elements is presented in the adjacency list form. In this representation, the first number in a line that is outside the brackets shows the network element under consideration. Numbers inside the brackets correspond to the network elements adjacent to the element of interest.

Groups 1-11 include different distribution substations, but may share generators. Specifically, groups 1 and 6 are connected by generator 68; groups 2, 6 and 10 by generator 69; groups 3, 10, and 11 by generator 79; groups 3 and 9 by generator 80;

groups 4, 7, and 11 by generator 54; groups 4, 10, and 11 by generator 55; groups 4 and 11 by generator 66; groups 5 and 7 by generators 70 and 83; groups 5 and 11 by generator 71; groups 5, 10, and 11 by generator 82; groups 8 and 10 by generator 64; groups 8 and 9 by generator 81; groups 9 and 11 by generator 58; and groups 10 and 11 by generators 74 and 77.

The following assumption is made: if two groups share a generator, this generator can supply power it produced to both groups, but cannot transport power from other generators to either of the groups. That is, a generator cannot operate as an interconnection.

Groups 1-5 cannot be simplified further as they already contain only one distribution substation (equivalently, one load). In group 6, both loads see generators 68 and 69 in the same way, that is, as a sub-topology with two VT-links adjacent to a VB-link. This sub-topology corresponds to the topology of group 3.

Groups 7-11 are more complex. Their analysis, however, can be conducted in a manner discussed in Section 3 in application to the network shown in Fig. 1a. In fact, the network in Fig. 1a is group 8. Results for other groups will be reported at the conference.

# 5 Conclusions

The current paper reports the initial results in mapping the Florida power grid with 53 distribution substations (loads) and 31 power plants (generators) onto a set of smaller sub-topologies with multiple generators and a single load connected to the grid by a single link. The preliminary analysis revealed that the system can be disintegrated into eleven groups (clusters) of interconnected generators and distribution substations. Some of these groups already contain sub-topologies of interest, that is, those with one load connected to the network of generators with a single link. Further analysis is required for groups with more complex topologies.

The preliminary analysis revealed a difference between microgrids such as, for example, ones in vehicles, and utility systems. Microgrids are pre-designed networks and therefore, exhibit symmetries in their structure. Utility systems, however, are not designed as a whole system, but rather evolve following economic development in a region. From the point of view of the survivability analysis, it makes microgrids easier to analyze with less computational costs. As a result, it is easier for a designer to make a decision on how to enhance their survivability due to their topology. In contrast, the analysis of utility systems is much more costly computationally and more demanding with regard to the effort required to conduct the analysis. Topological complexity makes it more difficult to predict the outcome of unexpected faults in such systems, identify weak links in them, and suggest improvements. At the same time, the performance of such systems affects larger population and they are more vulnerable to massive unpredicted damage than microgrids. That is, such systems are in greater need of advanced survivability analysis.

# Acknowledgments

We thank I. Abou Hamad and B. Israels for providing the adjacency matrix for the Florida power grid and a draft version of Fig. 4, based on a figure in Ref. [8].

This work was supported in part by U.S. National Science Foundation Grant No. DMR-1104829.